\documentclass[11pt]{article}

\usepackage{amsfonts,amssymb}

\newcommand{\Z}{\mathbb{Z}}
\newcommand{\R}{\mathbb{R}}

\def\Inn{ {\rm Inn}  }
\def\lt{{\triangleleft}}

\usepackage[dvips]{graphicx}
\usepackage{color}

\newtheorem{theorem}{Theorem}[section]

\newtheorem{lemma}[theorem]{Lemma}
\newtheorem{proposition}[theorem]{Proposition}
 
\newtheorem{definition}[theorem]{Definition}

\newtheorem{example}[theorem]{Example}

\newtheorem{remark}[theorem]{Remark}

\def\Inn{ {\rm Inn}  }
\def\lt{\,{\triangleleft}\,}

\begin{document}

\title{
Symmetric Extensions of Dihedral Quandles \\
and Triple Points of Non-orientable Surfaces }

\author{
J. Scott Carter\footnote{Supported in part by
NSF Grant DMS  \#0603926.}
\\ University of South Alabama
\and
Kanako Oshiro \footnote{Supported in part by Grant-in-Aid for JSPS Research Fellowships for Young Scientists.}\\ Hiroshima University
\and
Masahico Saito
\footnote{Supported in part by NSF Grant DMS
 \#0603876.}
\\ University of South Florida
}

\maketitle

\begin{abstract}

Quandles with involutions that satisfy certain conditions, called good involutions, 
can be used to color non-orientable surface-knots. 
We use 
subgroups 
of signed permutation matrices to construct 
non-trivial 
good involutions 
on extensions  of  
odd order dihedral 
quandles.

For the smallest example $\tilde{R}_3$ of order $6$ that is an extension
of the three-element dihedral quandle $R_3$, various symmetric quandle homology groups are 
computed, and applications to the minimal triple point number of surface-knots are given.
\end{abstract}

\section{Introduction}

In this paper, we construct an extension 
$\tilde{R}_m$ 
 of the dihedral 
quandle  
$R_m$ 
with a non-trivial good involution for each odd 
positive integer $m=2n+1$.
The extensions $\tilde{R}_m$ we construct are not involutory, and in particular, not isomorphic to 
dihedral quandles. 
As an application,
such an extension is used to study  the 
minimal triple point numbers of non-orientable 
surface-knots in thickened $3$-manifolds. 
Detailed definitions will be given in Section~\ref{prelimsec}. 

A quandle is a set 
with 
a binary operation that is self-distributive: $(a\lt b) \lt c = (a \lt c) \lt (b \lt c)$ and satisfies two other properties. The algebraic structure mimics the Reidemeister moves, and consequently quandles are a fundamental tool in knot theory. 
Quandle cohomology theories~\cite{CEGS,CJKLS} have been 
constructed, 
 and applied to knots
by 
using 
quandle elements
as colors 
 and cocycles as weights 
to define quandle cocycle invariants. 
The same construction was applied to surface-knots using triple points 
of the projections, and 
a 
variety of applications have been found. 

The original definitions of quandle colorings and quandle cocycle invariants are dependent upon 
orientations of the 
diagrams,
and in particular,  the invariants were defined at first only for orientable surface-knots.
A quandle is called {\it involutory} if $(x \lt y)\lt y=x$ holds  for any elements $x,y$ of the quandle.
An involutory quandle is also called a {\it kei} \cite{Taka}, and has 
the   
 property that colorings
are defined for unoriented knots and surfaces. 
To generalize 
involutory quandles  
and quandle cocycle invariants 
for unoriented diagrams
and 
non-orientable surfaces,
quandles with good involutions were defined~\cite{K} and studied~\cite{KO,O}. 
Quandles with good involutions are called {\it symmetric quandles}.
Constructions of symmetric quandles have depended mainly on 
computer calculations, or by hand for specific families of quandles, 
such as dihedral quandles. 
In particular, it was shown in \cite{KO} that any dihedral quandle of odd order
only has 
the 
trivial (identity map) good involution. 
In this paper, we give a construction of symmetric quandles 
via extensions of dihedral quandles. Specifically, 
we prove that for any odd order dihedral quandle, there is an extension 
with a non-trivial good involution  
that is connected and is not involutory.

The smallest of such extensions is given by a 
two-to-one 
quandle homomorphism 
$\tilde{R}_3 \rightarrow R_3$,  onto the three-element dihedral quandle.
Various homology groups of this quandle $\tilde{R}_3$ are computed, and specific 
non-trivial cycles and cocycles are presented.
Applications are given for the minimal triple point numbers of non-orientable surfaces.

The
 {\it  minimal triple point number} $t(K)$ 
 of a knotted or linked surface $F$ is defined  to be 
 the smallest number
of the triple point numbers among all the diagrams of the surface-link $F$, 
and denoted by $t(F)$. Quandle cocycle invariants \cite{CJKLS, CKS, CKS:book}  are used for studies of 
minimal triple point numbers of orientable surface-links; Satoh and Shima~\cite{SatShi}
determined the 
minimal triple point number  
of the 
2-twist-spun trefoil, and  Hatakenaka~\cite{Hata} gave 
a 
lower
bound 
of 
$6$ 
for the 
$2$-twist-spun 
 figure-eight knot. 
Kamada~\cite{Kamada} proved that, for any positive integer $N$, there is an orientable
 $2$-knot  (an embedded  sphere, a spherical surface-knot) 
 $K$ with $t(K)>N$. His argument (Alexander modules) does not 
 immediately apply to higher genus surfaces or non-orientable surfaces.

 Iwakiri~\cite{Iwakiri} used quandle cocycle invariants to provide 
a surface-knots $K$  
 with the triple point canceling number $\tau(K)$ 
 as large as you please. 
  The {\it triple point canceling number} 
  is the minimal number of $1$-handles 
  needed 
  to change a surface-knot 
 into 
 another 
 with a projection that has no triple point. 
 It was pointed out by Satoh that, since $\tau(K) \leq t(K)$, Iwakiri's  result implies 
 that  for any positive integer $N$, there is an orientable
 surface-knot $K$ with $t(K)>N$. Iwakiri's results can be applied to higher genus {\it orientable surfaces},
 but not to non-orientable surfaces.

In \cite{KO,O,Sat}, 
large minimal triple point numbers of 
(not necessarily orientable, 
and two-component) 
surface-links are realized. 
Their arguments, however, do not immediately apply to 
surface-knots.
In this paper, we give surface-knots in thickened $3$-manifolds with arbitrary 
large minimal triple point numbers.

The paper is organized as follows. 
A summary of definitions and known results necessary for this paper
are 
given 
in Section~\ref{prelimsec}.
Explicit constructions are given in Section~\ref{extsec}
to prove the existence of symmetric extensions
of  
odd order dihedral quandles. 
Symmetric quandle homology groups are computed for the smallest such example
in Section~\ref{triplesec}, 
and applications to the minimal triple points are presented.

\section{Preliminaries}\label{prelimsec}

In this section we give a summary of necessary definitions and set up  the 
notation.  

\subsection{Symmetric quandles} 

\begin{definition} {\bf \cite{Br88,Joyce,Matveev,Taka}}
{\rm 
A {\it quandle\/}, $X$, is a set with a binary operation
$(x,y) \mapsto x\lt y$ such that 

\noindent
(I. {\sc idempotency}) for any $x \in X$, $x \lt x=x$, 

\noindent
(II. {\sc right-invertibility}) for any $x,y \in X$,
 there is a unique $z \in X$, denoted by $x \lt \bar{y}$, such that $ x=z \lt y$, and 

\noindent 
(III. {\sc self-distributivity}) for any $x,y,z \in X$, we have
 $(x\lt y) \lt z=(x \lt z) \lt {(y \lt z)}$.}
\end{definition}
Typical examples 
are conjugations 
of groups $x \lt y=y^{-1}xy$. 
Any subset of a group closed under conjugation, thus, is a quandle.
In particular for the dihedral group $D_{2m}$ 
for any positive integer $m$, 
the subset $R_{m}$ of reflections forms a quandle by conjugation. In this paper, we will concentrate on the case 
$m=2n+1$ --- an odd integer. 
It is known that $D_{2m}$ has a presentation $\langle x,y: x^2=y^m=(xy)^2=1 \rangle$. 
The reflections and rotations of the regular $m$-gon are written as $xy^j$ and $y^j$
  for $j=0, \ldots, m-1$, respectively.
  Since  $(xy^j)^{-1} (xy^i )(xy^j)=xy^{2j-i}$, 
   the quandle
   $R_m$ can be identified with $\Z_m$ with the operation $i \lt j=2j-i \pmod{m}$.

Let $G$ be a group, $H$ a subgroup, $s:G\rightarrow G$ an automorphism such that 
$s(h)=h$ for each $h\in H$.  
Define a binary operation  on $G$ by 
$a\lt b=s(ab^{-1})b.$
Then this defines a quandle structure on $G$. 
This passes to a well-defined quandle structure on the right cosets
$G/H$ that is given by
$Ha \lt  Hb= Hs(ab^{-1})b$. 
In particular, if $\zeta \in Z(H)\cap H$ where $Z(H) = \{ \zeta \in G: \zeta h=h\zeta  \ {\mbox{\rm for all }} \  h\in H \}$, then
$Ha \lt Hb=Hab^{-1}\zeta b$ defines a quandle structure. Let us denote the resulting quandle by $(G,H,\zeta)$.
This construction is found in \cite{Joyce,Matveev}. 
For $G=D_{2m}$ with $H=\langle x \rangle$ and $\zeta=x$, one computes 
$Hy^i \lt Hy^j=Hy^i y^{-j} x y^j = Hxy^{2j-i}=Hy^{2j-i}$, so that we have
$R_m=(D_{2m}, H, x)$. 

Any element $c \in X$   
of a quandle $X$,    defines  a 
 function $S_c : X \to X$  by $(x)S_c = x \lt c$ and 
 that is a quandle automorphism (by axioms II and III) and that is 
 called   a {\it symmetry of $X$}.
The set of symmetries 
$\{S_c | c \in X\}$ forms the
 {\it inner automorphism group of $X$} 
that  is denoted by Inn($X$). 
If $\Inn(X)$ acts transitively on $X$, then $X$ is said to be {\it connected.}  

\begin{definition}  {\bf \cite{K}} 
{\rm An involution $\rho: X \rightarrow X$ defined on a quandle is a {\it good involution} if 
$x\lt {\rho (y) }= x\lt {\overline{y}}$ and $\rho(x\lt y)=\rho(x)\lt y$.
Such a pair $(X, \rho)$ is called a quandle
with a good involution or a {\it symmetric quandle}.
}\end{definition}

The associated group \cite{FR} of a quandle $X$ is
$G_X = \langle x \in X : x \lt 
y= y^{-1} x y  \rangle$. 
The associated group, $G_{(X,\rho)}$ of a symmetric quandle $(X, \rho)$ is defined~\cite{K, KO}
 by
$G_{(X, \rho)} = \langle x \in X : x \lt y= y^{-1} x y, \ 
 \rho(x)=x^{-1} \rangle$.    
The natural map $\mu: X \rightarrow G_{(X, \rho)} $ is the composition of the inclusion map
$X \rightarrow F(X)$ and the
projection map $F(X) \rightarrow G_{(X, \rho)} $,  where $F(X)$ is the free group on $X$. 
For a quandle $X$, an $X$-set \cite{FRS} 
is a set $Y$ equipped with a right action of the associated group
$G_X$. 
For a symmetric quandle $(X, \rho)$, an $(X, \rho)$-set is a set $Y$ 
equipped with a right action of the
associated group $G_{(X, \rho)} $.
We denote by $yg$ or by $y \cdot g$ the image of an element $y \in Y$  
under 
the action 
of  
$g \in G_{(X, \rho)} $. 
The   
 following 
three 
formulas hold:
$y \cdot (x_1 x_2) = (y \cdot x_1) \cdot x_2, \ $  
$y \cdot (x_1 \lt x_2) = y \cdot (x_2^{-1} x_1 x_2), \ $ and
$y \cdot (\rho( x_1 ) ) = y \cdot (x_1^{-1} ) $, 
for $x_1, x_2 \in X$ and $y\in Y$.

\subsection{Homology theories for symmetric quandles}

 A cohomology theory
 of quandles  was defined~\cite{CJKLS} as  a modification
of rack cohomology theory~\cite{FRS}.
In this section we review homology groups for symmetric quandles defined in \cite{K},
see also \cite{KO}.

\begin{sloppypar}
Let $Y$ be an $(X, \rho)$-set 
which may be empty. 
Let $C_{n}(X)_Y$ be the free abelian group generated by 
$(y, x_1,\ldots,x_n)$, where $y \in Y$ and $x_1, \ldots, x_n \in X$. 
For a positive integer $n$, 
let 
$C_0(X)_Y=\Z (Y)$, the free abelian group 
generated by  
$Y$, and 
set $C_n(X)_Y=0$ otherwise. 
(If $Y$ is empty, then  define $C_0(X) = 0$).  
Define the boundary 
homomorphism
\mbox{$\partial_n: C_{n}(X)_Y \longrightarrow C_{n-1}(X)_Y$}
by
\begin{eqnarray*}
\lefteqn{\partial_{n}(y, x_1,\ldots,x_n)
=\sum_{i=1}^{n}(-1)^i[(y, x_1,x_2,\dots,x_{i-1},\hat{x_i}, x_{i+1},\ldots,x_{n})}\\
&-&(y \cdot x_i, x_1\lt x_i,x_2\lt x_i,\ldots,x_{i-1}\lt x_i,\hat{x_i},  x_{i+1},\ldots, x_{n})]
\end{eqnarray*}
for $n\geq 1$ and $\partial_{n}=0$ for $n\leq 1$. Then
$C_{*}(X)_Y =\{C_{n}(X)_Y,\partial_n\}$ is a chain complex \cite{FRS}. 
Let
$D_{n}^Q(X)_Y$ be the subgroup 
of $C_{n}(X)_Y$ generated by 
$\cup_{i=1}^{n-1} \{ (y, x_1,\ldots, x_n )  \ | \ x_i=x_{i+1} \}$, and 
let 
$D_{n}^\rho (X)_Y$ be the subgroup 
of $C_{n}(X)_Y$ 
generated by $n$-tuples of  the form
$$ (y, x_1,\ldots, x_n ) + (y \cdot x_i , x_1\lt x_i ,\ldots, x_{i-1}\lt x_i , \rho(x_i), x_{i+1}, \ldots, x_n )$$     where $ y \in Y, \ \ $  $x_1, \ldots, x_n \in X, \ \ $ and $ i\in \{ 1, \ldots, n-1\}.$ 
Then $\{ D_{n}^Q(X)_Y, \partial_n \}$ and $\{ D_{n}^\rho (X)_Y, \partial_n \}$ are 
subcomplexes of $C_n$ \cite{KO}, and chain complexes 
$C_*^{R} (X)_Y$, $C_*^{Q} (X)_Y$, $C_*^{R, \rho} (X)_Y$, $C_*^{Q, \rho} (X)_Y$
are defined, respectively, from chain groups 
$C_n^{R} (X)_Y=C_n(X)_Y$, $C_n^{Q} (X)_Y=C_n (X)_Y/ D_n^Q (X)_Y$, 
$C_n^{R, \rho} (X)_Y=C_n (X)_Y/ D_n^\rho (X)_Y$, 
$C_n^{Q, \rho} (X)_Y = C_n (X)_Y / (D_n^Q  (X)_Y+ D_n^\rho (X)_Y)$. 
Their respective homology groups \cite{KO} are denoted by
$H_*^{R} (X)_Y$,  $H_*^{Q} (X)_Y$, $H_*^{R, \rho} (X)_Y$, 
and  $H_*^{Q, \rho} (X)_Y$, respectively. 
When $Y=\emptyset$, this subscript is dropped. 
Corresponding cohomology groups are defined as usual, as well as 
(co)homology groups with other coefficient groups.
\end{sloppypar}

An {\it extension} of a quandle $X$ is a surjective quandle homomorphism 
$f: E \rightarrow X$ such that for any element of $X$, the cardinality of the inverse image by $f$ is constant. 
We also say that $E$ is an extension of $X$.
In \cite{CENS}, an interpretation of quandle $2$-cocycles 
was given in terms of extensions of quandles, 
in a  manner similar to group extensions
by group $2$-cocycles. 
It is, therefore, a natural question to ask for a relation between 
symmetric quandle $2$-cocycles and  extensions of symmetric quandles.
Here we observe such an interpretation. 

Let $(X, \rho)$ be a symmetric quandle,
and $A$ be an abelian group, 
and $\phi: X^2 \rightarrow A$ be a symmetric quandle $2$-cocycle.
Specifically, $\phi$ satisfies 
\begin{itemize}
\item
$\phi(x_1,x_1)=0,$ 
\item
$\phi(x_1, x_2) - \phi(x_1, x_3) - \phi(x_1 \lt x_2, x_3) + \phi(x_1 \lt x_3, x_2 \lt x_3)=0$ for any 
$x_1, x_2, x_3 \in X,$
\item 
$\phi(x_1, x_2)+\phi(\rho(x_1), x_2)=0,$ and
\item $\phi(x_1, x_2)+\phi(x_1 \lt  x_2, \rho(x_2) )=0.$
\end{itemize}  
An extension of a quandle $X$ by a quandle $2$-cocycle $\phi$,
 denoted by $X \times _\phi A$, 
 was defined in \cite{CENS}
by $(x, a)\lt(y,b)=(x \triangleleft  y, a + \phi(x,y) )$.  
Define $\tilde{\rho} : X \times_\phi  A \rightarrow X \times_\phi A $ 
by $\tilde{\rho} (x, a)=(\rho(x), -a)$. 

\begin{proposition}
$(X \times_\phi A, \tilde{\rho})$ is a symmetric quandle.
\end{proposition}
This 
follows from 
direct calculations.

\subsection{Colorings of surface-knots by symmetric quandles}

 A knot diagram {\rm for a classical knot  $(n=1)$ or for a 
 surface-knot 
 $(n=2)$ is the image of a general position map 
 from a closed $n$-manifold (collection of circles or surfaces) 
 into $\R^{n+1}$ with crossing information indicated
 by breaking the under-arc or under-sheet (see \cite{CS:book} for details). 
 Let a 
 surface diagram $D$ of a surface-knot $F$ be given. We cut the diagram further into {\it semi-sheets} by considering the upper sheets 
 also 
 to  
 be broken along the double point arcs. Observe that in the local picture of a branch point there are two semi-sheets, at a double point there are $4$ semi-sheets, and at a triple point, there are $12$ semi-sheets.

Let $(X,\rho)$ 
denote a 
symmetric 
quandle, and  let $Y$ 
denote 
an $(X, \rho)$-set.  
The surface diagram $D$ 
has 
elements of $X$ assigned 
to the sheets and elements of $Y$ assigned to regions separated by 
the projection such that the following conditions are satisfied.
\begin{itemize}
\item (Quandle coloring rule on over-sheets) 
Suppose that two adjacent semi-sheets coming from an over-sheet of $D$ about a double curve are labeled by $x_1$ and $x_2$. If the normal orientations are coherent, then $x_1=x_2$, otherwise $x_1= \rho(x_2)$. 
\item (Quandle coloring rule on under-sheets) 
Suppose that two adjacent 
under-sheets
$e_1$ and $e_2$ 
are separated along a 
double curve 
and 
are labeled by $x_1$ and $x_2$. 
Suppose that one of the two semi-sheets
coming from the over sheet 
of $D$, say $e_3$, is labeled by $x_3$.  We assume that 
a local 
normal orientation of $e_3$ 
points 
from $e_1$ to $e_2$.  If the normal orientations of $e_1$ and $e_2$ are coherent, then $x_1\lt x_3 = x_2$, otherwise $x_1\lt x_3= \rho(x_2)$. 
\item (Region colors)
Suppose that two adjacent regions $r_1$ and $r_2$ which are separated by a semi-sheet, say $e$, are
labeled by $y_1$ and $y_2$, where $y_1, y_2 \in Y$. Suppose that the semi-sheet $e$ is labeled by $x$.
 If the normal orientation
of $e$ points from $r_1$ to $r_2$,  then $y_1 \cdot x = y_2$. 
\item
 An equivalence relation (of 
a local normal orientation assigned to each semi-sheet and a quandle element associated to this local orientation) is generated by the following rule 
({\it basic inversion}\/):  
Reverse 
the normal orientation of a semi-sheet and 
change 
the element $x$ assigned the sheet by $\rho(x)$. 
\end{itemize}
A {\it symmetric quandle coloring}, 
 or an {\it $(X, \rho)_Y$-coloring, }
 of a surface-knot diagram is such an equivalence class of symmetric quandle colorings. 
 See Fig.~\ref{inversion}.

 \begin{figure}[htb]
\begin{center}
\begin{minipage}{132pt}
\begin{picture}(132,70)

\qbezier(7,22)(7,22)(95,22)
\qbezier(7,22)(7,22)(43,56)
\qbezier(131,56)(43,56)(43,56)
\qbezier(131,56)(95,22)(95,22)

\put(66,20){\vector(0,-1){9}}
\multiput(66,38)(0,-4){5}{\line(0,-2){2}}
\put(74,10){$x$}

\end{picture}
\end{minipage}
\hspace{1cm}$=$\hspace{1cm}
\begin{minipage}{130pt}
\begin{picture}(130,70)

\qbezier(7,22)(7,22)(95,22)
\qbezier(7,22)(7,22)(43,56)
\qbezier(131,56)(43,56)(43,56)
\qbezier(131,56)(95,22)(95,22)

\put(69,38){\vector(0,1){26}}
\put(77,65){$\rho(x)$}

\end{picture}
\end{minipage}

\end{center}
\caption{A basic inversion}
\label{inversion}
\end{figure}
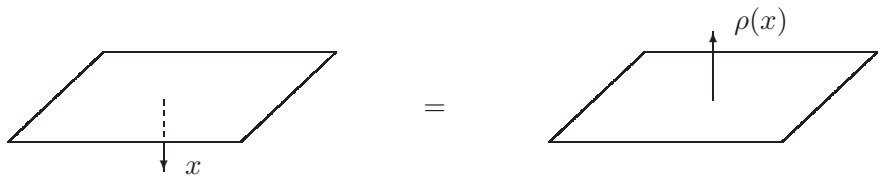

We call a diagram $D$ with an $(X,\rho)_Y$-coloring, $C_D$,  an {\it $(X,\rho)_Y$-colored diagram} and denote it by $(D,C_D)$.
Let $(D,C_D)$ and $(D',C_{D'})$ be $(X,\rho)_Y$-colored diagrams.

We say that $(D, C_D)$ and $(D', C_{D'})$ are {\it Roseman move equivalent} if they are related by a finite sequence of Roseman moves \cite{Rose} (see also \cite{CS:book})
such that the colors are kept constant outside of  each local move.

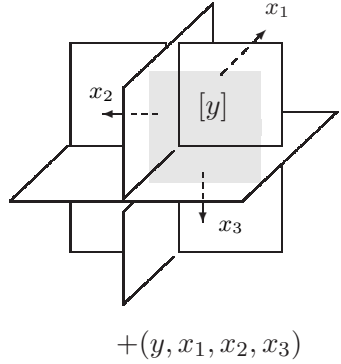
\begin{figure}[t]
\begin{center}
\begin{minipage}{130pt}
\begin{picture}(130,130)

\put(60,75){\setlength{\fboxrule}{18pt}\fcolorbox[gray]{.9}{.9}{}}
\thicklines
\thinlines

\qbezier(7,47)(7,47)(95,47)
\qbezier(7,47)(7,47)(29,68)
\qbezier(131,81)(95,47)(95,47)
\qbezier(131,81)(131,81)(110,81)

\qbezier(50,48)(50,48)(50,90)
\qbezier(84,122)(50,90)(50,90)
\qbezier(84,122)(84,122)(84,107)
\qbezier(50,48)(50,48)(69,66)

\qbezier(50,8)(50,8)(50,43)
\qbezier(50,8)(50,8)(69,26)
\qbezier(54,47)(50,43)(50,43)

\qbezier(71,68)(71,68)(110,68)
\qbezier(71,68)(71,68)(71,107)
\qbezier(71,107)(71,107)(110,107)
\qbezier(110,68)(110,68)(110,107)

\qbezier(71,28)(71,28)(110,28)
\qbezier(71,28)(71,28)(71,47)
\qbezier(110,28)(110,28)(110,61)

\qbezier(30,28)(30,28)(50,28)
\qbezier(30,28)(30,28)(30,47)

\qbezier(30,68)(30,68)(30,107)
\qbezier(30,107)(30,107)(66,107)
\qbezier(30,68)(30,68)(50,68)
\qbezier(66,106)(66,107)(66,107)

\put(34,86){\footnotesize\colorbox {white}{\textcolor{black}{$x_2$}} }

\put(50,80){\vector(-1,0){8}}
\multiput(55,80)(4,0){3}{\line(-1,0){2.0}}

\put(99,107){\vector(1,1){5}}

\qbezier(87,95)(87,95)(89,97)
\qbezier(91,99)(91,99)(93,101)
\qbezier(95,103)(95,103)(97,105)
\qbezier(99,107)(99,107)(101,109)

\put(80,47){\vector(0,-1){8}}
\multiput(80,50)(0,4){3}{\line(0,-1){2.0}}

\put(78,80){$[y]$}

\put(83,36){\colorbox [gray]{1.0}{\textcolor{black}{ \footnotesize$x_3$}}}
\put(100,118){ \footnotesize$x_1$}
\put(40,-10){ 
{ $+(y, x_1, x_2, x_3)$}}

\end{picture}
\end{minipage}

\end{center}
\caption{A weight of a triple point}
\label{weight}
\end{figure}

Let $(D,C_D)$ be an $(X,\rho)_Y$-colored diagram of an $(X,\rho)_Y$-colored surface-link $(F,C)$.
For a triple point of $D$, define the {\it weight} as follows:
Choose one of eight $3$-dimensional regions around the triple point and call the region a {\it specified region}.
There exist 12 semi-sheets around the triple points.
Let $S_T$, $S_M$ and $S_B$ be the three of them that face the specified region, where $S_T$, $S_M$ and $S_B$ are in the top sheet, the middle sheet and the bottom sheet at the triple point, respectively.
Let $n_T$, $n_M$ and $n_B$ be the normal orientations of $S_T$, $S_M$ and $S_B$ which 
point away from the specified region. 
Consider a representative of $C_D$ such that the normal orientations of $S_T$, $S_M$ and $S_B$ are $n_T$, $n_M$ and $n_B$ and
let $x_1$, $x_2$ and $x_3$ be the labels assigned the semi-sheets $S_B$, $S_M$ and $S_T$,  with the normal orientations  $n_B$, $n_M$ and $n_T$,  respectively. 
Let $y$ be the label assigned to the specified region. 
The weight is  
$\epsilon (y, x_1, x_2, x_3)$, where $\epsilon $ is $+1$ (or $-1$) if the triple of the normal orientations $(n_T,n_M,n_B)$ does (or does not, respectively) 
match 
the orientation of $\mathbb R^3$.
See Figure~\ref{weight}.
The sum $\sum_{\tau} \epsilon (y, x_1, x_2, x_3)$ of the weight over all 
triple points of a colored diagram $(D,C_D)$ represents a $3$-cycle $[c_D] \in C_3^{Q, \rho} (X)_Y$
\cite{KO}. 
Colored diagrams and cycles represented by colored diagrams are 
similarly defined 
when region colors are absent, or equivalently, when $Y=\emptyset$. 

\begin{lemma}{\bf \cite{KO,O} }\label{colorcyclelem}
If two colored diagrams $(D, C_D)$ and $(D', C_{D'})$ are 
Roseman  
move equivalent, then they represent homologous $3$-cycles,
$[C_D]=[C_{D'}]  \in H_3^{Q, \rho}(X)_Y$. 
\end{lemma}

For surface-knots or link in $\R^4$, it is known \cite{Rose} 
that two diagrams of equivalent (ambiently isotopic) surface-knot or a link 
are Roseman move equivalent.

\subsection{Triple point numbers}

Let $F$ be a surface-link and $D$ a diagram of $F$.
The minimal triple point number of $F$ is evaluated by quandle 
 invariants with symmetric quandle cocycles as follows:
 
\begin{lemma}{\bf \cite{KO,O}}\label{triplepointlem}
Let $(X,\rho)$ be a symmetric quandle.
Let $\theta:\mathbb Z( X^3)\to \mathbb Z$ be a symmetric quandle 3-cocycle 
$\theta \in C^3_{Q,\rho}(X)$
of $(X,\rho)$ such that $\theta (a,b,c)\in \{0,-1 ,1 \}$ for any $(a,b,c)\in X^3$.
For a symmetric quandle coloring $C_D$ of the diagram $D$, if $\theta ( [ C_D] ) =\alpha $ for $\alpha \in \mathbb Z$, then $t(F)\geq |\alpha |$. 
\end{lemma}
Especially, when a symmetric quandle $3$-cocycle $\theta :\mathbb Z(X^3) \to \mathbb Z$ of $(X,\rho)$ satisfies 
$\theta (a,b,c)=\pm 1$ for any $(a,b,c)\in X^3$ such that $a\not \in \{b, \rho(b) \}$ and $c\not \in  \{ b, \rho(b) \}$, we say $\theta $ is the $3$-cocycle with $\pm$ monic terms. 
 
For surface-knots and links  
in a thickened $3$-manifold
$M \times [0,1]$, where $M$ is a closed $3$-manifold, 
we take the natural 
 projection
 $p: M \times [0,1] \rightarrow M$ in the direction of the unit interval
 to define the diagrams. Then the minimal  triple point number is defined in the same manner 
 as above. 
By \cite{HN}, we may assume that diagrams in $M$ of 
equivalent (ambiently isotopic) surface-knot or link in $M \times [0,1]$ 
are Roseman move equivalent.

\section{Extensions of dihedral quandles with good involutions}\label{extsec}

In this section we prove:

\begin{theorem}\label{extthm}
For each positive integer $n$, there is an extension $\tilde{R}_{2n+1}$ 
of $R_{2n+1}$ with a non-trivial good involution $\rho$   
that is 
connected and 
is not involutory. 
\end{theorem}

Note that there are  extensions 
$R_{2(2n+1)} \rightarrow R_{2n+1}$,  and $R_{2(2n+1)}$ has a  non-trivial good involution for 
any $n>0$, see \cite{KO}.  In this case, however, $R_{2(2n+1)}$ is involutory and is not connected.
Connectedness of quandles play important roles in coloring knots and 
quandle homology. 
The identity map on an involutory quandle is a trivial good involution. 
Thus the interesting features of the extension $\tilde{R}_{2n+1}$ are its connectivity and its non-trivial good involution. 

The proof follows a  
construction of a group extension of the dihedral group and the definition of a quandle structure
on group cosets  
that were 
described 
in Section~\ref{prelimsec}.

\begin{definition}{\rm 
Let $e_j$ denote the  column vector  in $\R^m$  whose $j\/$th entry is $1$ and the remaining entries are  each $0$. A {\it signed permutation matrix} is a 
square matrix of size $m$ 
matrix  whose columns are 
of the form $(\pm e_{\sigma (1)}, \pm e_{\sigma (2)},\ldots, \pm e_{\sigma (m)})$ where $\sigma \in \Sigma_m$ is a permutation. The set of signed permutation matrices form a group $H_m$ of order $2^m m!$ that is called the {\it hyper-octahedral group}.  
Define the group $SH_m$ to be 
the signed permutation matrices of determinant $1$.
 }\end{definition}

 To avoid extra subscripts, we write $(\pm e_{\sigma (1)}, \pm e_{\sigma (2)},\ldots, \pm e_{\sigma (m)})$ as $(\pm {\sigma (1)}, \pm {\sigma (2)},\ldots, \pm {\sigma (m)}).$ 
 Then the matrix multiplication, in this notation, is written  as 
 $$(\epsilon_1 \cdot {\sigma (1)},\  \ldots,\  \epsilon_m \cdot {\sigma (m)})\cdot(\delta_1 \cdot {\tau (1)}, \ \ldots, \ \delta_m \cdot {\tau (m)})
=
(\epsilon_{\tau(1)} \delta_1\cdot \sigma (\tau(1)),  \ \ldots, \ \epsilon_{\tau(m)} \delta_m\cdot \sigma (\tau(m)))
,$$
where $\epsilon_i=\pm 1$ and $\delta_j=\pm 1$ for $i, j= 1, \ldots, m$. 
The product is obtained by looking at the entry in the $\tau(1)$ position of $(\epsilon_1 \cdot {\sigma (1)},\  \epsilon_2  \cdot {\sigma (2)},\ \ldots,\  \epsilon_m \cdot {\sigma (m)})$ and write that entry first after having been multiplied by $\delta_1$, then look at the entry in the $\tau (2)$ position and write that second after having been multiplied by $\delta_2$, and the process continues to the $m$th position.  
For example, 
$(1,5,4,-3,-2)\cdot(5,1,2,3,4)=(-2,1,5,4,-3)$ 
while 
$(5,1,2,3,4)\cdot(1,5,4,-3,-2)=(5,4,3,-2,-1).$

We identify a subgroup of $SH_m$ that maps onto the dihedral group.  
Let $m=2n+1$. Consider the subgroup $G_{2n+1}$   of $SH_{2n+1}$ that is generated by the pair of elements $a=(1, 2n+1, 2n, \ldots, n+2, -(n+1), \ldots, -3, -2)$ and $b=(2n+1,1,2,\ldots , 2n)$.

The dihedral group, 
$D_{2(2n+1)}$, will be identified  
with the image of its faithful representation 
in permutation matrices.
Specifically, 
we 
identify 
the reflection 
$x$ with $(1, 2n+1, 2n, \ldots, 2)$, 
the rotation 
$y$ with $(2n+1,1,2,\ldots , 2n)$, 
and $D_{2(2n+1)}$ with the subgroup of permutation matrices generated by 
these two elements. 
Then the map that takes each matrix $(a_{ij})$ to 
$(|a_{ij}|)$ defines  
a group homomorphism onto the dihedral group: 
$f: G_{2n+1} \rightarrow D_{2(2n+1)}$, 
such that $f(a)=x$ and $f(b)=y$.

\begin{lemma}\label{orderlem} The order of $G_{2n+1}$ is $(2n+1)\cdot 2^{2n+1}.$ 
The centralizer of $a$, $C(a)=\{ c \in G_{2n+1}: ac=ca\}$,  is generated by the elements 
$$(1,\;  \epsilon_2 \;(2n+1),\;  \epsilon_3 \ (2n), \; \ldots, \; \epsilon_{n+1} \ (n+2), \; - \epsilon_{n+1} \ (n+1), \; \ldots ,\;  -\epsilon_2\ (2))$$ where $\epsilon_j = \pm 1$ for  $j=2,\ldots , n+1.$ In particular, 
$|C(a)|= 2^{n+1}.$ 
\end{lemma}
{\it Proof.} 
Let $I_{\vec{\epsilon}}=(\epsilon_1 \, (1),\  \epsilon_2 \, (2),\  \ldots,\  \epsilon_{2n+1} \, (2n+1))$ where 
$\epsilon_j = \pm 1$ for $j=1, \ldots, 2n+1$,  
such that $\prod_{j=1}^{2n+1}\epsilon_j =1$ (an even number of entries are negative). There are 
$${{2n+1}\choose{0}} + {{2n+1}\choose{2}} + \cdots + {{2n+1}\choose{2n}} = 2^{2n}$$ such elements. 
We show that these elements are in $G_{2n+1}$. 
Since $a^2=(1, -2,-3, \ldots, -(2n+1)) $, 
for $i=1,\ldots ,2n$, 
$a^2 (b^{-i} a^2 b^i )$ has  
negative signs 
at the first and the 
$(i+1)$th entries, and 
positive signs otherwise. 
Hence $b^{-j} (a^2 b^{-i} a^2 b^i) b^j $ has 
negative signs 
at the 
 $(j+1)$th and 
 $(i+j+1)$th entries, and 
positive signs elsewhere.
By multiplying elements of these forms, any $I_{\vec{\epsilon}}$ with an 
even number of negative signs can be obtained. 
Since $b^ia=ab^{-i} I_{\vec{\epsilon}}$ for such an $ I_{\vec{\epsilon}}$, 
any element of $G_{2n+1}$ is written uniquely as $b^j  I_{\vec{\epsilon}}$ or $ab^j  I_{\vec{\epsilon}}$.
In Lemma~\ref{normal}, we will describe the multiplication of such normal forms. 
This is analogous to elements of  $D_{2n+1}$ having the form $y^j$ and $xy^j$. 
In total, we have $|G_{2n+1}|=(2n+1) \cdot 2^{2n+1}.$ 

If $ac=ca$, then $c$ must be of the form $I_{\vec{\epsilon}}$ or $a I_{\vec{\epsilon}}$. 
{}From the equation
$a\cdot a  I_{\vec{\epsilon}} = a I_{\vec{\epsilon}} \cdot a$, we have
$\epsilon_j = \epsilon_{2n+3-j}$. Since the determinants of the matrices 
are all $+1$, the initial sign, $\epsilon_1$ must be positive; or else, an odd number of the remaining $\epsilon_j$ are negative, but these signs agree in pairs. Thus 
$$a I_{\vec{\epsilon}} = 
(1, \; \epsilon_2 \ (2n+1), \; \ldots, \; \epsilon_{n+1} \ (n+1), \; -\epsilon_{n+1} (n) ,\; \ldots,\;  - \epsilon_2 \  (2)).$$
A similar computation gives that 
$I_{\vec{\epsilon}}=(1, \ \epsilon_2 \ (2), \; \ldots, \; \epsilon_{n+1} \ (n+1), \; \epsilon_{n+1} \ (n+2) ,\; \ldots, \; \epsilon_2 \ (2n+1))$ are the only diagonal signed permutation matrices that 
commutes with $a$.  
These are products of the $a I_{\vec{\epsilon}} \/$  that commute with $a$. This completes the proof. $\Box$

\bigskip

 Let $H=C(a)$ denote the centralizer of $a$. Consider the quandle structure $\tilde{R}_{2n+1}= (G_{2n+1}, H, a)$  given by $Hu \triangleleft  Hv = H uv^{-1}av$. 
 {}From the preceding lemma, we have $|\tilde{R}_{2n+1} |= (2n+1) 2^n $.

\begin{lemma} \label{qhomlem}
There is a surjective quandle homomorphism $f:\tilde{R}_{2n+1} \rightarrow R_{2n+1}$.  \end{lemma}
{\it Proof.} 
The group homomorphism $f: G_{2n+1} \rightarrow D_{2(n+1)} $ defined earlier,
by $f(a)=x$ and $f(b)=y$, 
satisfies $f(C(a))=\langle x \rangle$. 
Hence $f$ induces a quandle homomorphism (written by the same letter) 
$f: \tilde{R}_{2n+1}=(G_{2n+1}, C(a), a) \rightarrow R_{2n+1}=(D_{2(2n+1)}, \langle x \rangle, x)$. 
$\Box$

\begin{lemma} \label{goodinvlem}
For any positive integer $n$, the quandle $\tilde{R}_{2n+1}= (G_{2n+1}, H, a)$ 
has a good involution.
\end{lemma}
{\it Proof.\/} 
Define a map $\rho : \tilde {R}_{2n+1} \to \tilde {R}_{2n+1}$ by 
\[
\rho(Hu)=
\left\{
\begin{array}{ll}
H(-1,-2,\ldots, -(n+1),n+2,\ldots, 2n+1)u&\mbox{ if $n$ is an odd number,}\\
H(1,-2,\ldots, -(n+1),n+2,\ldots, 2n+1)u&\mbox{ if $n$ is an even number.}
\end{array}
\right.
\]
Note that the elements inserted 
$(-1,-2,\ldots, -(n+1),n+2,\ldots, 2n+1)$ and $(1,-2,\ldots, -(n+1),n+2,\ldots, 2n+1)$,
respectively, are indeed elements of $G_{2n+1}$, as they have even numbers of 
negative signs. 
We prove that this map is a good involution of $\tilde {R}_{2n+1}$.\\
(i) It is an involution by 
$$( \rho\circ \rho) (Hu)=H(\varepsilon 1,-2,\ldots, -(n+1),n+2,\ldots ,2n+1)^2u=Hu,$$
where $\varepsilon $ is the sign $\pm $. \\ 
(ii) For any element $Hu$ and $Hv$ in $\tilde {R}_{2n+1}$,
\begin{eqnarray*}
\rho(Hu)\triangleleft Hv&=& H(\varepsilon 1,-2,\ldots, -(n+1),n+2,\ldots ,2n+1)u\triangleleft Hv \\
 & = & H(\varepsilon 1,-2,\ldots, -(n+1),n+2,\ldots ,2n+1)uv^{-1}a v.
 \end{eqnarray*}
On the other hand,
$$\rho(Hu \triangleleft Hv)=\rho(Huv^{-1}a v)=H(\varepsilon 1,-2,\ldots, -(n+1),n+2,\ldots ,2n+1)uv^{-1}a v.$$
Hence, $\rho(Hu)\triangleleft Hv=\rho(Hu\triangleleft Hv)$ is satisfied.\\
(iii) For any element $Hu$ and $Hv$ in $\tilde {R}_{2n+1}$,
\begin{eqnarray*}
Hu\triangleleft \rho(Hv) & = & Hu\triangleleft H(\varepsilon 1,-2,\ldots, -(n+1),n+2,\ldots ,2n+1)v\\
 & = & Huv^{-1}(\varepsilon 1,-2,\ldots,  -(n+1),n+2,\ldots ,2n+1) \\
 &  & \quad a\  (\varepsilon 1,-2,\ldots,  -(n+1),n+2,\ldots ,2n+1)v \\
 & = & Huv^{-1}a^{-1}v.
 \end{eqnarray*}
The last equality is satisfied by 
$$(\varepsilon 1,-2,\ldots,  -(n+1),n+2,\cdots ,2n+1)\ a\  (\varepsilon 1,-2,\ldots,  -(n+1),n+2,\cdots ,2n+1)=a^{-1}.$$
The equality
$$( Hu\triangleleft \rho(Hv) ) \triangleleft Hv=Huv^{-1}a^{-1}v\triangleleft Hv=Huv^{-1}a^{-1}vv^{-1}av=Hu$$
implies $Hu\triangleleft \rho(Hv)=Hu\triangleleft \overline{Hv}.$
$\Box$

\begin{lemma} \label{noninvolem}
The quandle  $\tilde{R}_{2n+1} $ is not involutory
for any positive integer $n$.
\end{lemma}
{\it Proof.\/}
It is sufficient to show 
 that the equality $(Ha \lt Hb) \lt Hb = Ha$ does not hold in $\tilde{R}_{2n+1}$. 
Since 
$(Ha \lt Hb) \lt Hb=Hab^{-1}abb^{-1}ab=Hb^{-1}a^2b$, 
we  show  
that $b^{-1}a^2b \not \in H$.
One computes 
\begin{eqnarray*}
b^{-1}a^2 b &=& (2,\ldots,2n+1,1)(1,-2,-3,\ldots, -(2n+1))(2n+1,1,\ldots,2n)\\
            &=& (2,\ldots ,2n+1,1)(-(2n+1),1,-2,-3, \ldots , -2n)\\ 
            &=& (-1,2,-3, \ldots , -(2n+1)).
\end{eqnarray*}
By Lemma~\ref{orderlem}, $H$ is generated by some square matrices of size $2n+1$ whose $(1,1)$-entries are $1$, 
 and $(1,i)$-entries and $(j,1)$-entries are $0$ for any $i,j \not =1$.
Thus 
the matrix $(-1,2,-3, \ldots , -(2n+1))$ is not an element of $H$.
Therefore, $\tilde{R}_{2n+1}$ is not involutory.
$\Box$

\bigskip

Let ${\cal I}$ be the kernel of 
 the map $f: G_{2n+1}\to D_{2(2n+1)}$ 
 in Lemma~\ref{qhomlem}, 
 i.e., $${\cal I}=\{I_{\vec{\epsilon}}=(\epsilon_1(1),\ldots , \epsilon_{2n+1}(2n+1) ) ~|~\epsilon_i =\pm 1, ~\prod _{i=1}^{2n+1}\epsilon _i=1\}.$$
Define the maps $f_a: I \to I$, $f_b: I\to I$ and $f_{b^{-1}}: I \to I$ as follows:
\begin{eqnarray*}
f_a[(\epsilon_1(1),\ldots , \epsilon_{2n+1}(2n+1) )]
&=& (\epsilon_1(1),\epsilon_{2n+1}(2) , \ldots , \epsilon_{2}(2n+1) ), \\ 
f_b[(\epsilon_1(1),\ldots , \epsilon_{2n+1}(2n+1) )]
&=& (\epsilon_{2n+1}(1),\epsilon_1(2),\ldots , \epsilon_{2n}(2n+1) ),\  {\rm and} \\ 
f_{b^{-1}}[(\epsilon_1(1),\ldots , \epsilon_{2n+1}(2n+1) )]
&=& (\epsilon_{2}(1),\ldots ,\epsilon_{2n+1}(2n), \epsilon_{1}(2n+1) ), 
\end{eqnarray*} 
where the order of $(\epsilon_2, \ldots, \epsilon_{2n+1} )$ is reversed for $f_a$, 
and $(\epsilon_1, \ldots, \epsilon_{2n+1} )$ is cyclically permuted for $f_b$ and $f_{b^{-1}}$. 
We can easily see that $f_a$, $f_b$ and $f_{b^{-1}}$ are  
automorphisms of ${\cal I}$, 
$f_a^{-1}=f_a$ and $f_b^{-1}=f_{b^{-1}}$. Moreover, their actions on the diagonal matrices correspond to the dihedral actions of the elements $x$, $y$, and $y^{-1}$. To imagine this action consider a necklace of $2n+1$ pearls an even number of which are black. The automorphisms $f_b$ and $f_b^{-1}$ act as rotations while $f_a$ acts a reflection of the necklace that fixes the first pearl on the strand.

We also consider distinguished elements  
\begin{eqnarray*}
I_{+} &=& (1,\ldots,n,-(n+1),n+2,\ldots ,2n, -(2n+1)), \ {\rm and} \\ 
I_{-} &=& (-1, 2,\ldots, n+1, -(n+2), n+3, \ldots ,2n+1), 
\end{eqnarray*} 
which have exactly two minus signs.  
Observe 
that $f_b(I_+)=I_-$. 
Then the following equalities hold: 
$${I_{\vec{\epsilon }}}\; a=af_a(I_{\vec{\epsilon }}), \quad 
{I_{\vec{\epsilon }}}\, b=bf_b(I_{\vec{\epsilon }}), \quad 
 {I_{\vec{\epsilon }}}\, b^{-1}=b^{-1}f_b^{-1}(I_{\vec{\epsilon }}), $$
$$ \ ba=ab^{-1}I_+, \quad
{\rm and}  \quad  \ b^{-1}a=abI_-. $$
Consequently, we have the following relations:
$$b^{\pm j} a =a b^{\mp j} \left[ {\prod}_{k=0}^{j-1} f^{\mp k}_b (I_{\pm}) \right], 
\quad  {\rm and} \quad   
I_{\vec{\epsilon}} \, b^{\pm j} = b^{\pm j} f_b^{\pm j} (I_{\vec{\epsilon}}).$$
If $1 \le j \le n$, then  \vspace{2.5mm} 
${\prod}_{k=0}^{j-1} f^{-k}_b (I_{+})$
is a diagonal matrix 
in $\cal{I}$ that has 
 two blocks of $j$ contiguous $(-)$-signs; the first block ends at  $n+1$, 
and the second block ends at  $2n+1$.  
\vspace{2.5mm} 
In particular, ${\prod}_{k=0}^{n-1} f^{ -k}_b (I_{+})$ 
has exactly one $(+)$-sign at position $(1)$.  \vspace{2.5mm} 
Similarly,
${\prod}_{k=0}^{n-1} f^{ k}_b (I_{-})$ 
has exactly one $(+)$-sign at position $(n+1)$. 

\vspace{0.25cm}  
 
For $i=1, \ldots, 2n+1$, let $I(i)$ denote the diagonal matrix that has 
 exactly one $(+)$-sign at position $(i,i)$ and $(-1)\/$s at the other positions along the diagonal.

\begin{lemma}\label{normal}
The following product formulas hold:
\begin{enumerate}
\item  \hfil $(b^i I_{\vec{\epsilon}})(b^j I_{\vec{\delta}}) = b^{i+j} f_b^j(I_{\vec{\epsilon}}) I_{\vec{\delta}},$ \hfill \hfill
\item   \hfil $(b^i I_{\vec{\epsilon}})(a b^j I_{\vec{\delta}}) = ab^{j-i}
\left[ {\prod}_{k=0}^{i-1} f^{-k}_b (I_{+}) \right] 
 f_b^j(f_a(I_{\vec{\epsilon}})) I_{\vec{\delta}},$ \hfill \hfill
 \item   \hfil $(a b^i I_{\vec{\epsilon}})(a b^j I_{\vec{\delta}}) = a^2 b^{j-i} \left[ {\prod}_{k=0}^{i-1} f^{j-k}_b (I_{+}) \right] f_b^j(f_a(I_{\vec{\epsilon}})) I_{\vec{\delta}}, $ \hfill \hfill  
 \item   \hfil $(a b^i I_{\vec{\epsilon}})( b^j I_{\vec{\delta}}) = a b^{i+j}
f_b^j(I_{\vec{\epsilon}}) I_{\vec{\delta}}.$ \hfill \hfill
\end{enumerate}
\end{lemma}
{\it Proof.} The calculations follow directly from the formulas above. $\Box$

\begin{lemma}\label{connlem} 
The quandle $\tilde{R}_{2n+1}$ is connected.
\end{lemma}
{\it Proof.\/}
Since any element of $G_{2n+1}$ is written as $b^{i}I_{\vec {\epsilon }}$ or $ab^{i}I_{\vec {\epsilon }}$ for some $i\in \{0,\ldots ,2n\}$ and 
$I_{\vec{\epsilon}}=(\epsilon_1(1),\ldots, \epsilon_{2n+1}(2n+1)) \in {\cal I},$ 
any element of $\tilde{R}_{2n+1}$ is written as $Hb^{i}I_{\vec {\epsilon }}$.  
We further abbreviate $I_{\vec{\epsilon}}$ as $\vec{\epsilon}$. Thus Greek letters with arrows in the formulas below indicate diagonal matrices. 

\bigskip

\noindent 
\underline{Claim~1} :   {\it 
 For any $Hb^{i}
 {\vec {\epsilon }} \in \tilde{R}_{2n+1}$, there exists a matrix $
 {\vec{\delta }}\in \cal{I}$ such that 
$$ ( Hb^i\lt Hb^{n+i}) \lt Hb^{n+i} 
{\vec{\delta}}=Hb^{i}
{\vec {\epsilon }}.$$ 
 }

\noindent 
{\it Proof of Claim~1}.
Using Lemma~\ref{normal}, one computes 
\begin{eqnarray*}
\lefteqn{
(H b^i\lt H b^{n+i} ) \lt H b^{n+i} 
{\vec{\delta}} } \\
& =  & H b^i b^{-(n+i)}ab^{n+i}  \, 
{\vec{\delta}} \, b^{-(n+i)} a b^{n+i} \, 
{\vec{\delta}}\\
&=& 
H b^{-n}
a b^{n+i} \, {\vec{\delta}}  \, b^{-n-i} a b^{n+i} \,  {\vec{\delta}}  \\  
& = & 
H a b^{n} \left[\prod_{l=0}^{n-1}  f_b^{l}(I_-)\right] 
 b^{n+i} \, {\vec{\delta}}  \, b^{-n-i} a b^{n+i} \,  {\vec{\delta}}  \\  
& = & H  b^{n} I(n+1)  
 b^{n+i} \, {\vec{\delta}}  \, b^{-n-i} a b^{n+i} \,  {\vec{\delta}}  \\ 
& = & H  b^{2n+i} f_b^{n+i}(I(n+1))  
 \, {\vec{\delta}}  \, b^{-n-i} a b^{n+i} \,  {\vec{\delta}}  \\ 
& = & H  b^{n} f_b^{-n-i}\left ( f_b^{n+i}(I(n+1))  
 \, {\vec{\delta}}\, \right)  a b^{n+i} \,  {\vec{\delta}}  \\
& = & H  b^{n} I(n+1)  
 f_b^{-n-i} (\, {\vec{\delta}} \,)  a b^{n+i} \,  {\vec{\delta}}  \\
 & = & H  b^{n} a b^{n+i} (f_b^{n+i} \circ f_a)\left( I(n+1)  
 f_b^{-n-i} ( \,{\vec{\delta}}\,) \right)  {\vec{\delta}} \\
  & = & H a b^{-n}   \left[\prod_{l=0}^{n-1}  f_b^{-l}(I_+)\right]  b^{n+i} f_b^{n+i}( I(n+2))  
 (f_b^{n+i} \circ f_a \circ f_b^{-n-i}) ( \,{\vec{\delta}}\,)  {\vec{\delta}} 
 \\
  & = & H  b^{i}  f_b^{n+i}(I(1)) (f_b^{n+i}( I(n+2))  
 (f_b^{n+i} \circ f_a \circ f_b^{-n-i}) ( \,{\vec{\delta}}\,)   {\vec{\delta}} \\
   & = & H  b^{i}  I(n+i+1) I(i+1)  
 (f_b^{n+i} \circ f_a \circ f_b^{-n-i}) ( \,{\vec{\delta}}\,)   {\vec{\delta}} \\
    & = & H  b^{i}  I(n+i+1) I(i+1)  
 (f_b^{2n+2i} \circ f_a ) ( \,{\vec{\delta}}\,)   {\vec{\delta}} \, .
 \end{eqnarray*}
 Then one computes 
$$
( f_{b}^{2n+2i} \circ f_a ) (\, {\vec{\delta}} \, ) =(\delta_{2i} , \delta_{2i-1}, \ldots ,\delta_1, \delta_{2n+1}, \ldots, \delta_{2i+1} ), $$ 
where we abbreviated the notation omitting the numbers that specify the entry, as 
all matrices in question for the rest of the proof are diagonal
(for example, the first entry of the first case should read $\delta_{2i} (1)$.
With this convention, we obtain the expression
\begin{eqnarray*}
C(\, {\vec{\delta}}\,) 
&=&
 I(n+i+1) I(i+1)    ( f_{b}^{2n+2i}  \circ f_a ) (\, {\vec{\delta}} \, ) {\vec{\delta} } \\  
&=& 
 I(n+i+1) I(i+1)  (\delta_1 \delta_{2i},  \delta_2 \delta_{2i-1 }, 
\ldots , \delta_{2i } \delta_1, \delta_{2i+1} \delta_{2n+1}, \ldots,  \delta_j \delta_k, \ldots 
  \delta_{2n+1} \delta_{2i+1} )  
\end{eqnarray*}
where the generic  $k\/$th term  $ \delta_j \delta_k$ has the property that $j+k \equiv 2i+1 \pmod{2n+1}$.  For  any choice of ${\vec{\delta}}$, 
the entry 
of  $ ( f_{b}^{2n+2i}  \circ f_a ) (\, {\vec{\delta}}\,) \cdot ( \, {\vec{\delta} }\,)$ 
 for which $2j \equiv 2i+1 \pmod{2n+1}$
 is $(\delta_j)^2$.  
 Solving, we see that this is position $j=i+n+1$. Therefore, 
\begin{eqnarray*}
\lefteqn{
C(\, \vec{\delta} \,)  =  (
 \delta_{1} \delta_{2i} , 
\delta_{2} \delta_{2i-1},  \ldots , \delta_{i} \delta_{i+1}, -\delta_{i+1} \delta_{i}, \delta_{i+2} \delta_{i-1} 
, \ldots, 
\delta_{2i} \delta_1,  }\\   
& & 
\delta_{2i+1} \delta_{2n+1} , 
\ldots, 
 \delta_{n+i} \delta_{i+n+2}, 
  - (\delta^2_{n+i+1}), 
  \delta_{n+i+2}\delta_{i+n}, 
   \ldots, \delta_{2i-1} \delta_{2n+1} ).  \end{eqnarray*}

\medskip

We choose a set of coset representatives $\{ H b^i \vec{\epsilon}\,  \} $
for each $i=1, \ldots, 2n+1$  as follows.
Let ${\vec{\epsilon}} \in {\cal{I}}$ be given. Define $S_{\vec{\epsilon}}\,(i)  \subset \{i+2, \ldots, n+i+1\}$ 
by the condition $s\in S_{\vec{\epsilon}}\,(i)$ if and only if $\epsilon_s = \epsilon_{2i+2-s}$ where throughout all subscripts
 are taken mod $(2n+1)$ but  
chosen to be the representative element in $\{1, \ldots, 2n+1 \}$. 
We show that 
$H b^i {\vec{\epsilon}} $ and $H  b^i {\vec{\eta}} $ represent the same coset if and only if 
 $S_{\vec{\epsilon}}\,(i) = S_{\vec{\eta}}\,(i).$
The cosets  $H b^i {\vec{\epsilon}} $ and $H  b^i {\vec{\eta}} $ are the same if and only if
$b^i  {\vec{\epsilon}}\,  \cdot {\vec{\eta}}\ b^{-i} \in H$,
and $b^i  {\vec{\epsilon}}\,  \cdot {\vec{\eta}}\ b^{-i} =
f^{-i}_b ( \, {\vec{\epsilon}}\,  \cdot {\vec{\eta}}\,) $.
On the other hand, from the proof of Lemma~\ref{orderlem}, 
a diagonal matrix is an element of $H$ if and only if 
it has the form 
$(1, \ \epsilon_2' ,  \ldots,  \epsilon_{n+1}' , \epsilon_{n+1}' , \ldots, \epsilon_2' ) $, so that 
 $f^{-i}_b( \, {\vec{\epsilon}}\,  \cdot {\vec{\eta}}  \,)\in H$ if and only if  
 $\epsilon_{i+1} \cdot \eta _{i+1}=1$, $\epsilon_{i+2} \cdot \eta _{i+2} =\epsilon_i \cdot \eta _{i}$, and so  
 forth, 
 which implies  $S_{\vec{\epsilon}\, \cdot \vec{\eta}} \, = \{i+2, \ldots, n+i+1\}$. 
 Thus $\epsilon_{i+1} = \eta_{i+1}$, and 
 $S_{\vec{\epsilon}}\,(i) = S_{\vec{\eta}}\,(i).$ 
Hence  $\{ S_{\vec{\epsilon}}\,(i) \} $ represent   cosets uniquely.

Now we show that for any ${\vec{\epsilon}} \, \in {\cal{I}}$ there is $C(\, \vec{\delta} \,) $
such that $H b^i {\vec{\epsilon}}$ and $H b^i C(\, \vec{\delta} \,) $ represent the same coset.
If  $n+i+1 \in S_{\vec{\epsilon}}\,(i)$, then $\epsilon_{n+i+1}=\epsilon_{n+i+2}$.
Hence  
the product $ \delta_{n+i+2} \delta_{i+n} = - (\delta_{n+i+1})^2$ is negative. 
This  sign  then determines the sign of the entry in position $(i+n)$ of $C(\delta)$. 
We continue in this way:  $i+n \in S_{\vec{\epsilon}}(i)$ if and only if $\delta_{n+i+3}\delta_{i+n-1} = -1$.
Inductively, the 
signs of the products $\delta_k \delta_{2i+1-k}$ are
 determined 
(cyclically) 
to the right of the $(n+i+1)\/$st entry 
by the values to the left and  by considering whether or not a given element is in $S_{\vec{\epsilon}}\,(i)$.
In this way, we 
obtain families 
$\vec{\delta}$  
for which 
$S_{C(\vec{\delta})}\,(i) = S_{\vec{\epsilon}}\,(i)$.

\bigskip

\noindent 
\underline{Claim~2} :  {\it For any elements $Hb^i 
{\vec{\epsilon_1}}\, $ and $Hb^i 
{\vec{\epsilon_2}}\, $ of $\tilde{R}_{2n+1}$, there exist  symmetries $S_1,\ldots , S_\mu $ of $\tilde{R}_{2n+1}$
 such that $(Hb^i 
 {\vec{\epsilon_1}}) 
 (S_1^{e_1}\circ \cdots \circ S_\mu^{e_\mu} ) 
= Hb^i 
{\vec{\epsilon_2}}$, 
where $e_j=\pm 1 $ for $j=1, \ldots, \mu$. } 

\bigskip

\noindent 
{\it Proof of Claim~2}. 
{}From Claim~1,  it follows that for each $j=1,2$, 
there is a matrix $
{\vec{\delta_j}}\, $ such that 
$(Hb^i\lt Hb^{n+i})\lt Hb^{n+i} 
{\vec{\delta_j}}\, =Hb^{i} 
{\vec {\epsilon _j}}\, $. 
Denote the symmetries coming from $Hb^{n+i}$, $Hb^{n+i} 
{\vec{\delta_1}}\, $ and $Hb^{n+i} 
{\vec{\delta_2}}\,$ by $S_b$, $S_1$ and $S_2$, respectively. Then 
\[
((Hb^i 
{\vec{\epsilon_1}}\, )S_1^{-1})S_2 =((((Hb^i 
{\vec{\epsilon_1}})S_1^{-1})S_b^{-1})S_b)S_2=((H b^i )S_b)S_2= Hb^i   
{\vec{\epsilon_2}}
\]
as desired. 

\bigskip

\noindent
\underline{Claim~3}:
{\it 
For any element $Hb^i  
{\vec{\epsilon}}\, $ of $\tilde{R}_{2n+1}$, there exists some symmetries $S_1,\ldots ,S_\nu  $  of $\tilde {R}_{2n+1}$ such that 
$(H)(S_1^{e_1} \circ \cdots \circ S_\nu^{e_\nu} ) =Hb^i 
{\vec{\epsilon}}$,
where $e_j=\pm 1$ for $j=1, \ldots, \nu$. 
 } \\

\noindent
{\it Proof of Claim~3}. 
Reduce the integer $k = i(n+1)$ modulo $2n+1$. 
Then 
we have
$$
H\lt H b^k= Hb^{-k} a b^k = Hb^{2k}f_b^k\left[ \prod_{l=0}^{k-1}f_b^{l}(I_-) \right] =Hb^i f_b^k\left[ \prod_{l=0}^{k-1}f_b^{l}(I_-)\right].
$$
To simplify the notation, let 
 ${\vec{\alpha}}\, = 
 f_b^k(\prod_{l=0}^{k-1}f_b^{l}(I_y)).$  
By Claim~1, there exists some symmetries $S_1,\ldots ,S_\mu $ of $\tilde {R}_{2n+1}$ 
such that $(H b^i 
{\vec{\alpha}}\, ) (S_1^{e_1} \circ \cdots \circ S_\mu^{e_\mu}) =H b^i 
{\vec{\epsilon}}\,$,
where $e_j=\pm 1$ for $j=1,  \ldots , \mu$. 
Therefore 
we obtain 
\[
(H) 
(S_{H b^k}  \circ S_1^{e_1} \circ \cdots \circ S_\mu^{e_\mu} ) 
 =H b^i 
{\vec{\epsilon}}\,
\]
as desired. 

\bigskip

\noindent 
Lemma follows from Claim~3. $\Box$

\bigskip

\noindent 
Theorem~\ref{extthm} follows from Lemmas~\ref{qhomlem},
 \ref{goodinvlem},  \ref{noninvolem}, 
and \ref{connlem}.

\begin{example} {\rm
The extension $\tilde {R}_3=(G_3, C(a), a)$, where $H=C(a)$, 
$a=(1,3,-2)$
and $b=(3,1,2)$, 
 consists of $6$ elements.
The six elements are represented by 
$0$ through $5$ as
$(0=H,\ 1=Hb^2, \ 2=Hb, \ 3=H(-1,-2,3), \ 4=Hb^2(-1,2,-3), \ 5= Hb(-1,-2,3) )$
 with the quandle operation given by the following  table.
\begin{center}
\begin{tabular}{c|cccccc} 
$R\triangleleft  C$ & $0$ & $1$ & $2$ & $3$ & $4$ & $5$ \\ \hline
$0$ & $0$ & $5$ & $1$ & $0$ & $2$ & $4$ \\
$1$ & $2$ & $1$ & $3$ & $5$ & $1$ & $0$ \\
$2$ & $4$ & $0$ & $2$ & $1$ & $3$ & $2$ \\
$3$ & $3$ & $2$ & $4$ & $3$ & $5$ & $1$ \\
$4$ & $5$ & $4$ & $0$ & $2$ & $4$ & $3$ \\
$5$ & $1$ & $3$ & $5$ & $4$ & $0$ & $5$ 
\end{tabular}
\end{center}
The map $f: \tilde {R}_3 \rightarrow R_3$ is given by $f(i)\equiv i \pmod{3}$. 
The good involution is the involution $\rho=(0~3)(1~4)(2~5)$. 
}\end{example}

\section{Homology groups of $\tilde{R}_3$ and triple point numbers} \label{triplesec}

For $\tilde{R}_3$, computer calculations give the results below
on symmetric quandle homology groups. 
Let $\chi_{(x,y,z)} \in C^3_{Q, \rho}(\tilde{R}_{2n+1}, \Z)$ be the characteristic function. 
Define a $3$-cochain 
\begin{eqnarray*}
A(x,y,z)  &=& 
\chi_{(x,y,z)} 
- \chi_{ (\rho(x),y,z)} 
- \chi_{(x\lt y ,\rho(y),z)} 
- \chi_{ (x \lt z, y \lt z, \rho(z) )} \\
 & &  
 +\chi_{ (\rho(x) \lt y ,\rho(y),z)} 
 + \chi_{ (\rho(x)\lt y , y \lt z, \rho(z))} \\ 
 & & 
 + \chi_{  ((x \lt y) \lt z,  \rho(y)\lt z, \rho(z))} 
-\chi_{ ( (\rho(x) \lt y) \lt z,  \rho(y) \lt z, \rho(z))} .
\end{eqnarray*}
Then {\it Mathematica} calculations show:

\begin{lemma} \label{keylemma}
Let $\tilde{R}_{3}$ be as above. 
\begin{itemize}
\setlength{\itemsep}{-3pt}
\item[{\rm (i)}]
 $H_2^{Q, \rho}(\tilde{R}_{3}, \Z)=0$,  $H_3^{Q, \rho}(\tilde{R}_{3}, \Z)\cong \Z$. 
\item[{\rm (ii)}]
The $3$-chain  $c=(2,1,2)+
(2,0,1) 
-(1,0,2)-(0,2,1) \in C_3^{Q, \rho}(\tilde{R}_{2n+1}, \Z)$
is a $3$-cycle $( c \in  Z_3^{Q, \rho}(\tilde{R}_{3}, \Z) )$ that represents a generator
$[c]$ of  $H_3^{Q, \rho}(\tilde{R}_{3}, \Z)\cong\Z$. 
\item[{\rm (iii)}]
Any   $3$-cycle with less than $4$ basis terms (triples) is null-homologous.
\item[{\rm (iv)}]
The $3$-cochain  $\phi= A(0,1,0)+A(0,1,2)-A(0,2,1)$
is a $3$-cocycle $( \phi \in Z^3_{Q, \rho}(\tilde{R}_{3}, \Z) )$ that represents a generator
of  $H^3_{Q, \rho}(\tilde{R}_{3}, \Z)\cong \Z$ dual to $[c]$, 
that is: 
$\phi([c])=1$.
\item[{\rm (v)}]
The $3$-cochain $\phi'=A(0,1,0)+A(0,1,2)+A(0,2,0)-A(0,2,1)
+A(1,0,1)-A(1,0,2)+A(1,2,0)+A(1,2,1)
+A(2,0,1)+A(2,0,2)-A(2,1,0)+A(2,1,2) $  is a $3$-cocycle 
with $\pm$ monic terms such that $\phi' ([c])=4$.
\end{itemize}
\end{lemma}

\begin{figure}[htb]
\begin{center}
\includegraphics[width=4in]{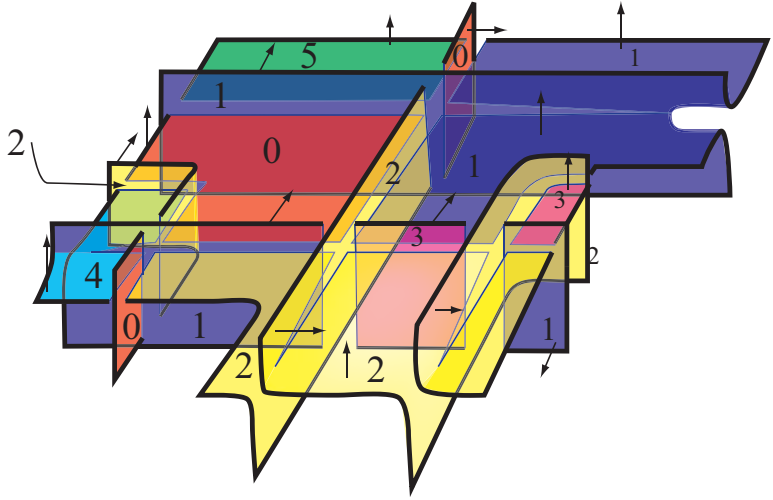} 
\vspace{-7mm}
\end{center}
\caption{A diagram of the surface constructed} 
\label{sfcediag}
\end{figure}

\begin{theorem}\label{TPthm}
 For any positive integer $N$, 
 there is a closed $3$-manifold $M$ and a non-orientable 
  surface-knot $F$ in $M \times [0,1]$ such that 
 $t(F) >N$.
 \end{theorem}
 {\it Proof.\/} 
 For 
 the $3$-cycle $c$ of Lemma~\ref{keylemma} (ii), 
make a colored triple point in a $3$-ball for each basis term. 
The degenerating terms are capped by branch points.
Connect them together to form a larger $3$-ball $B$ with all triple points and branch points included, see Fig.~\ref{sfcediag}.
The boundary $\partial B$ contains a colored classical link diagram illustrated in 
Fig.~\ref{doublehandle2}.
One can also obtain Fig.~\ref{doublehandle2} from 
the formula 
for the $3$-cycle $c$ of Lemma~\ref{keylemma} (ii)
 as follows:
 The $3$-cycle $c$ is represented by a colored diagram with region colors 
 as depicted in  Fig.~\ref{3cycle}. 
Take the ``double'' of Fig.~\ref{3cycle} and extend, see Fig.~\ref{double}.
By smoothing the black dots 
that represent branch points, 
we obtain Fig.~\ref{doublehandle2}.

\begin{figure}[htb]
\begin{center}
\includegraphics[width=3in]{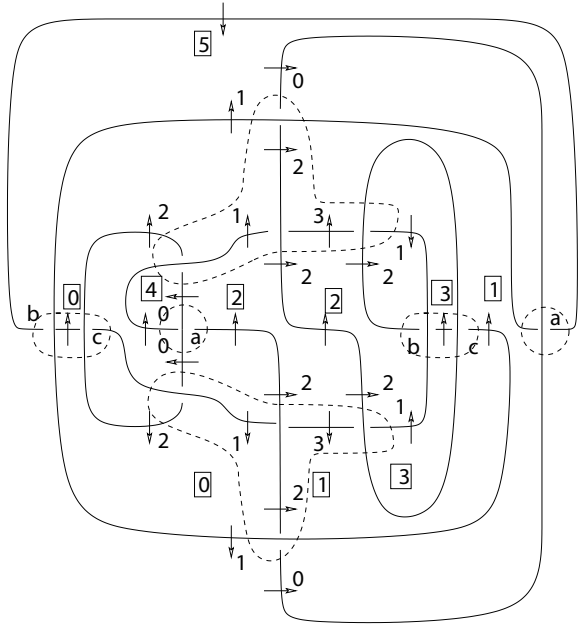}
\end{center}
\caption{Adding $1$-handles} 
\label{doublehandle2}
\end{figure}

Then add $1$-handles to connect double curves of the diagram. In Fig.~\ref{doublehandle2}
the attaching disks of $1$-handles are indicated by dotted circles. 
The shapes of the circles, T-shaped, oval and circle, respectively, together with 
the colors of arcs indicate the pairs of the attaching regions. Note that the oval and 
circle ones must be rotated $180$ degrees before identifying.
This twist   
makes the surface non-orientable. 
A handlebody $H$  of genus $3$ results 
as an ambient manifold, and it has $5$ closed curves on the boundary.

\begin{figure}[htb]
\begin{center}
\includegraphics[width=2in]{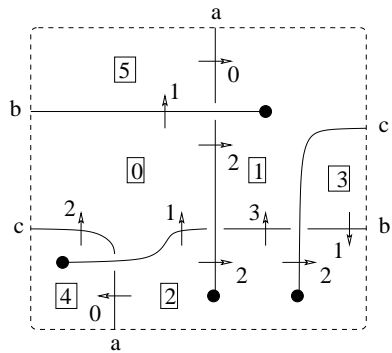}
\end{center}
\caption{Representing the $3$-cycle $c$}
\label{3cycle}
\end{figure}

\begin{figure}[htb]
\begin{center}
\includegraphics[width=2.5in]{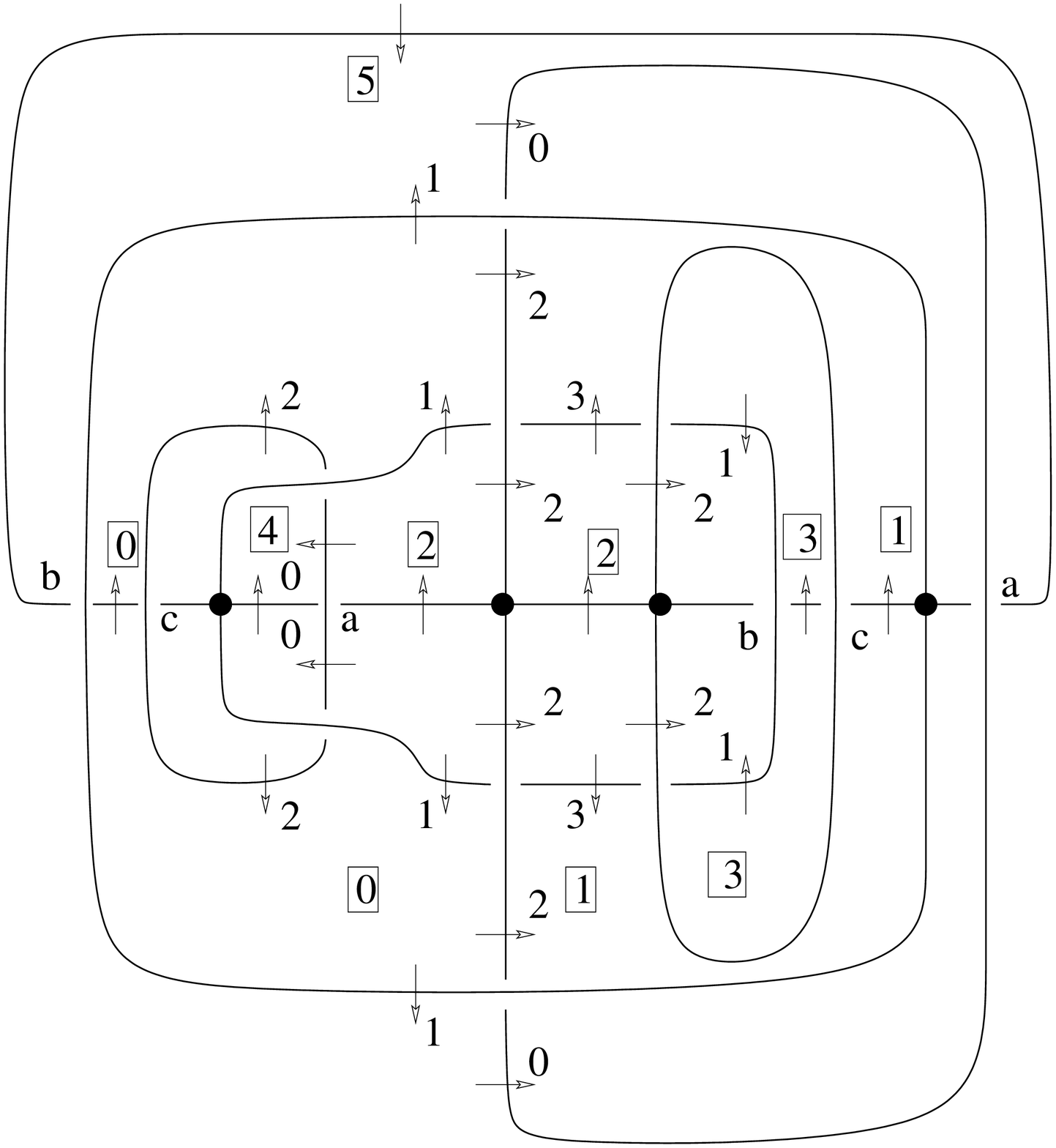}
\end{center}
\caption{Assembling triple points} 
\label{double}
\end{figure}

Attach 
$2$-handles to $H$ along the closed curves on the boundary.
Let $M'_0$ be the result, which is a compact $3$-manifold 
with boundary.
By capping off the boundary of $M'_0$ by handlebodies, we obtain a closed orientable
$3$-manifold $M_0$ with a diagram $D_0$ in it, that represents $c$.
Let  $m$ be an integer such that $4m>N$. 
Taking 
an $m$-fold 
knot connected sum, we have a connected closed $3$-manifold 
$M=\#_m M_0$  
and a connected, colored diagram $D=\#_m D_0$ 
 in $M$ which represents $m c$.
By lifting $D$ 
 to $M \times [0,1]$, we obtain the surface-knot $F$ whose minimal triple point number is greater than $N$:
Using the 
$3$-cocycle 
$\phi '$ in Lemma~\ref{keylemma},  
we have 
$t(F)\geq 4m$
by Lemma~\ref{triplepointlem}. 
Therefore we obtain the inequality $t(F)> N$.
 $\Box$

Note that using the $3$-cocycle $\phi '$, we can also prove that the minimal triple point number of the constructed surface-knot 
$F$ 
in the above proof is exactly $4m$.

The next result shows that homological conditions on the surface 
changes the triple point numbers. 

\begin{proposition}\label{8prop}
Any surface-knot diagram colored with $\tilde{R}_3$ in 
a closed $3$-manifold $M$  that is null-homologous 
in $H_2(M;\Z_2)$ as an underlying generic surface,  
and with less than $8$ triple points, is null-homologous in $H_3^{Q, \rho}(\tilde{R}_3, \Z)$.
\end{proposition}

For the proof, we need the following lemma, calculated by {\it Mathematica}.
Let $Y=\{  \alpha, \beta \}$, and let $\tilde{R}_3$ act on $Y$ by $\alpha\cdot u=\beta$, $\beta\cdot u =\alpha$
 for any $u \in \tilde{R}_3$.
This provides cycles represented by colored diagrams with regions with checkerboard colorings.
The map of deleting the first factor 
$\pi:  (\alpha {\mbox{ \rm or}}\ \beta, x_1, \ldots, x_n) \mapsto (x_1, \ldots, x_n)$
induces a chain map $\pi: C_n^{Q, \rho}(\tilde{R}_3, \Z)_Y \rightarrow C_n^{Q, \rho}(\tilde{R}_3, \Z)$.

\begin{lemma} \label{checkerlemma}
Let $\tilde{R}_3$, $Y$ be as above. 
\begin{itemize}
\setlength{\itemsep}{-3pt}
\item[{\rm (i)}]
 $H_3^{Q, \rho}(\tilde{R}_3, \Z)_Y \cong \Z \times \Z_3 $. 
\item[{\rm (ii)}]
The $3$-chain  
\begin{eqnarray*}
\gamma &=& 
(\alpha, 0,1,0) + (\alpha, 0,1,2 ) + (\alpha, 0,2,0 ) + (\alpha, 1,2,0) \\
& & - (\alpha, 2,1,0) + (\beta, 0,2,0)  + (\beta, 1,2,0) + (\beta, 2,0,1) 
 \in C_3^{Q, \rho}(\tilde{R}_3, \Z)_Y
 \end{eqnarray*}
is a $3$-cycle $( \gamma \in  Z_3^{Q, \rho}(\tilde{R}_3, \Z)_Y)$ that represents a generator
$[\gamma ]$ of  $\Z \subset H_3^{Q, \rho}(\tilde{R}_3, \Z)_Y$. 
\item[{\rm (iii)}]
Any  $3$-cycle with  less than  $8$ basis terms (quadruples) is null-homologous. 
\item[{\rm (iv)}]
The $3$-cochain $\phi''=A(0,1,0)+A(0,1,2)+A(0,2,0)-A(0,2,1)+A(1,0,1)-A(1,0,2)+A(1,2,0)+A(1,2,1)+A(2,0,1)+A(2,0,2)-A(2,1,0)+A(2,1,2)$ is a $3$-cocycle 
$(\phi'' \in  Z_3^{Q, \rho}(\tilde{R}_3, \Z) )$
with $\pm$ monic terms such that $\phi'' \circ \pi_*([\gamma ])=8$. 
\end{itemize}
\end{lemma}

\begin{lemma} \label{twicelemma}
The induced map 
$\pi_*: H_3^{Q, \rho}(\tilde{R}_3, \Z)_Y \rightarrow H_3^{Q, \rho}(\tilde{R}_3, \Z)$
restricted 
to the $\Z$ factor multiplies the generator by $2$. 
\end{lemma}
{\it Proof.\/}  One computes 
$\phi \circ \pi_*( [\gamma]) = \phi ( (0,1,0) + (0,1,2 ) + (0,2,0 ) + (1,2,0)
- (2,1,0) + ( 0,2,0)  + (1,2,0) + (2,0,1) )=2$. 
$\Box$

\bigskip

\noindent
{\it Proof (of Proposition~\ref{8prop}).}
Let $D$ be a colored diagram in 
a closed $3$-manifold $M$ whose underlying generic surface represents a null-homologous class in 
$H_2(M;\Z_2)$, and 
that is non-trivial in $H_3^{Q, \rho}(\tilde{R}_3)$. 
In particular, we have $\phi(D)\neq 0$, where $\phi$ is a cocycle in Lemma~\ref{keylemma} (iv).
Then there is a checkerboard coloring for $D$ 
as it is null homologous in $H_2(M;\Z_2)$, 
and
let $D'$ be the cycle 
in $Z_3^{Q, \rho}(\tilde{R}_3, \Z)_Y$ represented by $D$ with the checkerboard coloring.
Since $\pi_*( [D'] )\neq 0$ in $H_3^{Q, \rho}(\tilde{R}_3, \Z)$,
by Lemma~\ref{twicelemma}, $[D']$ is non-trivial in  $H_3^{Q, \rho}(\tilde{R}_3, \Z)_Y$.
Then Lemma~\ref{checkerlemma} (iii) implies that $D$ must have at least $8$ triple points.
$\Box$
 
\bigskip

\begin{figure}[htb]
\begin{center}
\includegraphics[width=2in]{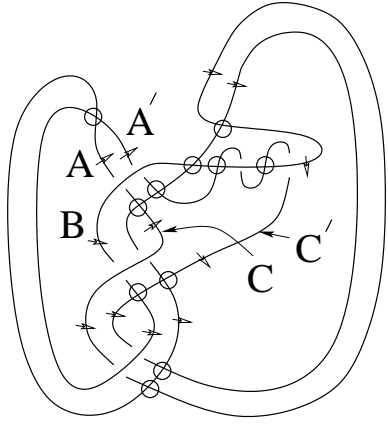}
\end{center}
\caption{A $3$-colorable virtual knot that has no non-trivial coloring by $\tilde{R}_3$} 
\label{virt}
\end{figure}

\begin{remark} {\rm 
Note that any coloring by $\tilde{R}_3$ gives rise to a coloring by $R_3$ 
by the map $f: \tilde{R}_3 \rightarrow R_3$, but the converse is not necessarily true.
All classical $3$-colorable knots we tested, however, 
have non-trivial colorings by $\tilde{R}_3$,
so we conjecture that it is always the case.

On the other hand, there are virtual knots that are $3$-colorable but are not 
non-trivially colored by $\tilde{R}_3$. Such an example is depicted in Fig.~\ref{virt}.
The virtual knot in the figure is $3$-colorable, and any coloring by $R_3$ is 
determined uniquely by the colors on the arcs labeled $A$ and $B$,
so that there are $9$ colorings by $R_3$, three of which are trivial. 
Suppose there is a non-trivial coloring by $\tilde{R}_3$.
If the induced $3$-coloring is trivial, say $0 \in R_3$, 
then the coloring consists of 
the  
two 
lifted 
colors, say $0$ and $3$.
These two elements, however, satisfy $0\lt 3=0$ and $3 \lt 0=3$, so that a connected virtual 
knot will be monochromatic, a contradiction. Hence we may assume that the given 
non-trivial coloring induces a non-trivial $3$-coloring. 
Let $\alpha$, $\alpha '$ and $\beta \in \tilde{R}_3$
be the colors assigned on the arcs $A$, $A'$ and $B$, respectively, 
with respect to the right direction nomals 
as depicted. The induced 
colors of $R_3$ are the same for $\alpha$ and $\alpha'$, 
that is, $f(\alpha)=f(\alpha')$. Hence $\alpha'=\alpha$ or $\rho(\alpha)$.
Note, by inspection, that $ (x \lt y) \lt y = \rho(x)$ holds for any $x, y \in \tilde{R}_3$
such that $x \neq y$ and $x \neq \rho(y)$. 
Hence the colors of the arcs on $C$ and $C'$ are 
$\alpha \lt \beta$ and $( (\alpha' \lt \beta )\lt  \beta ) \lt \beta  
=\rho(\alpha' \lt \beta )=\rho(\alpha') \lt \beta$, respectively. By tracing this arc back we see that for any choice of $\alpha$ or $\alpha'$, no consistent coloring can be obtained.
} \end{remark}

\noindent 
{\bf Concluding remarks.}
The most remarkable aspect of this quandle $X=\tilde{R}_3$  
   is its free part in $H_3$
despite 
its being 
connected 
(Lemma~\ref{keylemma} (i)). 
It is known \cite{LN} that the ordinary quandle homology groups do not have free part 
if it is connected. This shows that the symmetric quandle homology is 
quite different from the original quandle homology, and this fact should be useful 
for non-orientable surfaces. 
We  conjecture, however,  that any surface-knot diagram in $\R^3$ colored by $\tilde{R}_3$
represents null-homologous class in $H_3^{Q, \rho}(\tilde{R}_3, \Z)_Y$.
It is an interesting fact that, from Proposition~\ref{8prop}, the homology class  
a  surface represents
in homology groups of the parent $3$-manifold
is related to the non-triviality in quandle homology and the minimal triple point number.

It is an interesting problem to compute the quandle (co)-homology of quandle extensions (which are given by surjective  quandle homomorphisms) in terms of the homological information of the source, target, and fiber.

\subsection*{Acknowledgments} We are grateful to Professors Seiichi Kamada and Shin Satoh for 
numerous valuable comments and information.

\end{document}